\documentclass[11pt]{article}
\usepackage{amsmath,amssymb,amsthm}

\def\qq{q^{-1}}
\def\ot{\otimes}
\def\vv{V^{\otimes 2}}
\def\K{\mathbb{K}}

\def\al{{\alpha}}

\def\S{\mathcal{S}}

\def\RRR{{\mathcal{R}}}
\def\RR{R^{-1}}

\def\be{\begin{equation}}
\def\ee{\end{equation}}

\begin{document}

\title{On Symmetrizers in Quantum Matrix Algebras}
\author{\rule{0pt}{7mm} Dmitry Gurevich\thanks{gurevich@ihes.fr}\\
{\small\it Institute for Information Transmission Problems}\\
{\small\it Bolshoy Karetny per. 19,  Moscow 127051, Russian Federation}\\
\rule{0pt}{7mm} Pavel Saponov\thanks{Pavel.Saponov@ihep.ru}\\
{\small\it
National Research University Higher School of Economics,}\\
{\small\it 20 Myasnitskaya Ulitsa, Moscow 101000, Russian Federation}\\
{\small \it and}\\
{\small \it
Institute for High Energy Physics, NRC "Kurchatov Institute"}\\
{\small \it Protvino 142281, Russian Federation}\\
\rule{0pt}{7mm} Vladimir Sokolov\thanks{sokolov@itp.ac.ru}\\
{\small\it Institute for Information Transmission Problems}\\
{\small\it Bolshoy Karetny per. 19,  Moscow 127051, Russian Federation}}

\maketitle

\noindent
{\bf 1.} In this note we are dealing with a particular class of quadratic algebras, which are called Quantum Matrix algebras. Let us recall that by a quadratic algebra
one  means a quotient algebra  $A=T(V)/\langle J \rangle$, where $V$ is a vector space ($\dim V=N<\infty $), $T(V)=\bigoplus_{k\geq 0} V^{\ot k}$ is the free tensor algebra,
generated by $V$, and $\langle J \rangle$ is the ideal generated  by a subspace $J\subset \vv$. Algebraic properties of such algebras and of their inhomogeneous counterparts
are intensively discussed in mathematical literature (see for instance \cite{PP}). In this note we consider the problem of constructing some projectors in homogenous components of
such algebras, which are analogs of the usual symmetrizers. In the classical case such symmetrizers are defined via the  flip  $P$. Below, the role of $P$ is played by the so-called 
Hecke symmetries. Let us recall that an operator $R:\vv\to\vv$ is called a braiding if it is subject  to the braid relation:
$$
R_{12} \,R_{23}\,R_{12}=\,R_{23}\, R_{12}\, R_{23},\quad R_{12}=R\ot I,\,\, R_{23}=I \ot R.
$$
Hereafter, $I$ is the identity operator in $V$ (or another vector space) and  $q\in \K$ is a parameter. The ground field $\K$ is assumed to be of characteristic 0.

If a braiding $R$  meets a  complementary relation $(R-q\, I)( R+\qq\, I)=0$, then $R$ is called a Hecke symmetry. The best known Hecke symmetries are those coming
from the quantum groups $U_q(sl(N))$ or their super-analogs but there exist Hecke symmetries, which are deformations neither of the flips nor of the super-flips (see  \cite{G}).
Also, in this paper there  were constructed symmetrizers for the $R$-symmetric and $R$-skew-symmetric algebras of the space $V$.

\medskip

\noindent
{\bf 2.} Other quadratic algebras, related to braidings $R$, are the RTT algebras and Reflection Equation (RE) algebras. The first of them is defined by the system
\be
R\,T_1\, T_2  =T_1\,T_2\, R,  \qquad T_1=T\ot I,\quad T_2=I\ot T,
\label{RTT}
\ee
where $T=\|t_i^j\|_{1\leq i,j \leq N}$ is an $N\times N$ matrix, its entries generate the RTT algebra.

By using  a second braiding $F$, in a sense compatible with $R$, it is also possible to define more general Quantum Matrix algebras (see \cite{IOP}).

The main objective of this note is to present a method, which hopefully enables one to construct symmetrizers in the homogenous components of the RTT algebras and to illustrate
it with two low-dimensional examples. A way of extending this method onto other Quantum Matrix algebras is discussed at the end of the note.

Note that the defining system (\ref{RTT}) can be written in the form:
$$
T_1\, T_2=\RRR(T_1\, T_2),\quad {\rm where}\quad \RRR(T_1\, T_2)=\RR\, T_1\, T_2\, R.
$$
If a braiding $R$ is a Hecke symmetry, the operator $\RRR:W^{\ot 2}\to W^{\ot 2}$, where $W=\mathrm{span}(t_i^j)$, has three eigenvalues $1$, $-q^2$ and $-q^{-2}$. Hence the operator
$$
\S=2_q^{-2} (\RRR+q^2\, I)(\RRR+q^{-2}\, I),\quad {\rm where}\quad k_q\stackrel{\mathrm{def}} {=}\frac{q^k-q^{-k}}{q-\qq}, \quad \forall\,k\in\mathbb{Z},
$$
is an idempotent, called  {\em symmetrizer}. It maps the space $W^{\ot 2}$ onto its subspace of  symmetric elements, i.e. such that $\S(z)=z$. Thus,  each element $z$ of the second homogenous component $A^{(2)}$ of the  quadratic algebra $A=T(W)/\langle J \rangle$ equals its  symmetrized form $\S(z)$, modulo the ideal $\langle J \rangle$,
where $J=\mathrm{Im} (\RRR-I)$.

In \cite{GPS} an analogous symmetrizer $\S^{(3)}:W^{\ot 3}\to W^{\ot 3}$  was constructed in the third homogenous component $A^{(3)}$ of the algebra $A$. Its explicit form is
\be
\S^{(3)}:=\al \S_{12}\,\S_{23}\, \S_{12}\,\S_{23}\, \S_{12}+\beta \S_{12}\,\S_{23}\, \S_{12}+\gamma\S_{12}=
\al\S_{23}\, \S_{12}\,\S_{23}\, \S_{12}\,\S_{23}+\beta \S_{23}\, \S_{12}\,\S_{23}+\gamma\S_{23},
\label{coeff}
\ee
with some coefficients $\al$, $\beta$, $\gamma$, such that $\al+\beta+\gamma=1$. Hereafter, lower indices indicate positions where  the operators are located.

Since the operator defined by formula (\ref{coeff}) has the property:
$$
\S^{(3)}=\S_{12}\,\S^{(3)}=\S^{(3)}\, \S_{12}= \S_{23}\,\S^{(3)}=\S^{(3)}\, \S_{23},
$$
it can be written in the form $\S^{(3)}=p(\S_{12}\, \S_{23})=p(\S_{23}\, \S_{12})$, where $p$ is a  polynomial of degree 3.

We have found the minimal polynomial of the operator $\S_{12}\, \S_{23}$. It equals
\be
m_3(x)=x(x-1)\Big(x-\frac{1}{2_q^2}\Big)\Big(x-\frac{(2_q^2-2)^2}{2_q^{4}}\Big).\label{tree}
\ee
Consequently, the operator $\frac{1}{\kappa}\, p(\S_{12}\, \S_{23})$, where $p(x)=m(x)(x-1)^{-1}$ and $\kappa=p(1)$, is the projector of the space $W^{\ot 3}$ onto its
subspace, corresponding to the eigenvalue 1. Note that formula (\ref{tree}) enables one to compute all coefficients in (\ref{coeff}).

\medskip

\noindent
{\bf 3.} Our next objective is to construct the symmetrizer $\S^{(4)}:W^{\ot 4}\to W^{\ot 4}$. Similarly to the previous case, we compute the minimal polynomial $m_4(x)$ for the operator
$$
\S^{(3)}_{123}\, \S^{(3)}_{243}:\,\,W^{\ot 4}\to W^{\ot 4}.
$$
It is $m_4(x)=x(x-1)(x-\nu_1)(x-\nu_2) ( x-\nu_3)$, where
$$
\nu_1=\frac{1}{3_q^2},\quad \nu_2=\frac{(2_q^2-2)^2}{4\cdot 3_q^2},\quad \nu_3=\frac{(3_q^2-3)^2}{4\cdot 3_q^4}.
$$
Then, substituting the product $\S^{(3)}_{123}\, \S^{(3)}_{234}$ into the polynomial $p_4(x)=m_4(x)(x-1)^{-1}$, we get (after a proper normalization) the symmetrizer
$S^{(4)}$, we are looking for. It was verified that by interchanging the factors in this product, we get the same result.  

Hopefully, an analogous procedure enables one to find the higher symmetrizers.
By assuming $\S^{(n)}$ to be known, we define $\S^{(n+1)}$ as follows
$$
S^{(n+1)}_{1,2,\dots,n+1}=\frac{1}{\kappa_{n}} \,p_n\Big(S^{(n)}_{1,2,\dots,n} S^{(n)}_{2,3,\dots,n+1}\Big),\qquad p_n(x)=\frac{m_n(x)}{ (x-1)},
$$
where $m_n(x)$ is the minimal polynomial of the operator $S^{(n)}_{1,2,\dots,n} S^{(n)}_{2,3,\dots,n+1}$ and $\kappa_n=p_n(1)$ is a normalizing factor. This procedure is based
on the conjecture that $1$ is a simple root of the polynomial $m_n(x)$ for any $n$. Presumably the degree of the polynomial $m_n(x)$ is $n+1$ for a generic $q$, but it becomes equal to 3 at the limit $q=1$.

This conjecture  is also justified by our computation in the case $n=5$. In this case  the roots of the polynomial $m_5(x)$  are $0,\,\,1$ and the following ones
$$
\nu_1=\frac{1}{4_q^2}, \quad \nu_2=\frac{(2\cdot 2_q^2-5)^2}{9\cdot 2_q^4}, \quad \nu_3=\frac{(2\cdot 2_q^2-5)^2}{9\cdot 4_q^2}, \quad
\nu_4=\frac{(4_q^2-4)^2}{9\cdot 4_q^4}.
$$
It would be interesting to find a regular pattern in the series of the polynomials $m_n(x)$.

Note that the existence of the symmetrizer $S^{(n)}$ entails the following claim. If a  Hecke symmetry $R=R(q)$ is an analytical matrix-function
in a vicinity of $q=1$, then in some vicinity of $q=1$ the dimensions of the homogenous components $A^{(n)}$ are equals to these corresponding to $q=1$.

In conclusion, we want to emphasize that all constructions above are valid for other Quantum Matrix algebras, in particular, for the RE algebras. In the latter case
we have only to replace the products $T_1\, T_2$ etc by $L_1\, R\, L_1\, R^{-1}$ etc. For detail the reader is referred to \cite{GPS}. In a similar manner we proceed for other
Quantum Matrix algebras, but the shift of the matrix $L_1$ to the higher positions should be done by means of the second braiding $F$ instead of $R$ (see \cite{IOP}).

\end{document}